\documentclass[submission]{FPSAC2024}
\pdfoutput=1 %% I HAVE ADDED THIS!!!

\usepackage{subeqnarray,bm}

\newtheorem{theorem}{Theorem}[section]
\newtheorem{proposition}[theorem]{Proposition}

\newcommand{\be}{\begin{equation}}
\newcommand{\ee}{\end{equation}}
\newcommand{\textbfit}[1]{\textbf{\textit{#1}}}
\def\reff#1{(\protect\ref{#1})}
\def\proof{\par\medskip\noindent{\sc Proof.\ }}

\def\sketchofproofof#1{\par\medskip\noindent{\sc Sketch of proof of #1.\ }}

\renewcommand{\qed}{ $\square$ \bigskip}

\newcommand{\eee}{{\rm ee}}
\newcommand{\ooo}{{\rm oo}}

\newcommand{\eqdef}{\stackrel{\rm def}{=}}

\newcommand{\scrd}{{\mathcal{D}}}

\newcommand{\scro}{{\mathcal{O}}}
\newcommand{\scrr}{{\mathcal{R}}}

\newcommand{\bfa}{{\mathbf{a}}}
\newcommand{\bfb}{{\mathbf{b}}}
\newcommand{\bfc}{{\mathbf{c}}}
\newcommand{\bfd}{{\mathbf{d}}}
\newcommand{\bfe}{{\mathbf{e}}}
\newcommand{\bfff}{{\mathbf{f}}}

\newcommand{\bfx}{{\mathbf{x}}}
\newcommand{\bfy}{{\mathbf{y}}}

\newcommand{\Z}{{\mathbb Z}}

\newcommand{\R}{{\mathbb R}}

\newcommand{\Sym}{{\mathfrak{S}}}

\newcommand{\myge}{\succeq}

\newcommand{\cyc}{{\rm cyc}}

\newcommand{\cpeakeven}{{\rm cpeakeven}}
\newcommand{\cpeakodd}{{\rm cpeakodd}}

\newcommand{\cdrise}{{\rm cdrise}}

\newcommand{\cdfall}{{\rm cdfall}}

\newcommand{\fix}{{\rm fix}}

\newcommand{\ba}{{\bm{a}}}

\newcommand{\bll}{{\bm{\ell}}}
\newcommand{\bp}{{\bm{p}}}
\newcommand{\bq}{{\bm{q}}}

\newcommand{\bt}{{\bm{t}}}

\newcommand{\bzero}{{\bm{0}}}

\def\Xhat{{\widehat{X}}}

\newcommand{\sn}{{\rm sn}}
\newcommand{\cn}{{\rm cn}}
\newcommand{\dn}{{\rm dn}}

\newcommand{\Sn}{{\rm Sn}}
\newcommand{\Cn}{{\rm Cn}}
\newcommand{\Dn}{{\rm Dn}}
\newcommand{\En}{{\rm En}}

\newenvironment{scarray}{
             \textfont0=\scriptfont0
             \scriptfont0=\scriptscriptfont0
             \textfont1=\scriptfont1
             \scriptfont1=\scriptscriptfont1
             \textfont2=\scriptfont2
             \scriptfont2=\scriptscriptfont2
             \textfont3=\scriptfont3
             \scriptfont3=\scriptscriptfont3
           
           \begin{array}{c}}{\end{array}}

\title[Coefficientwise Hankel-total positivity of the Schett polynomials]{Coefficientwise Hankel-total positivity \\ of the Schett polynomials}

\author{Bishal Deb\thanks{\href{mailto:bishal@gonitsora.com}{bishal@gonitsora.com}. Partially supported by the DIMERS project ANR-18-CE40-0033.}\addressmark{1,2}
\and
Alan D.~Sokal\thanks{\href{mailto:sokal@nyu.edu}{sokal@nyu.edu}}\addressmark{1,3}
}

\address{\addressmark{1}University College London, London, UK \\
         \addressmark{2}Sorbonne Universit\'e and Universit\'e Paris Cit\'e, CNRS, 
	 Laboratoire de Probabilit\'es, Statistique et Mod\'elisation, 
	 Paris, France\\
         \addressmark{3}New York University, New York, NY, USA
}

\received{\today {\bf NEED TO PUT DATE!!!}}

\abstract{We prove the coefficientwise Hankel-total positivity
of the even and odd subsequences of Schett polynomials $X_n(x,y,z)$.}


\usepackage[backend=bibtex]{biblatex}
\begin{document}

\maketitle
%% note that you DO NOT have to put your abstract here -- it is generated by \maketitle and the \abstract and \resume commands above

\section{Introduction}

A half-century ago, Schett \cite{Schett_76,Schett_77} introduced implicitly,
and then Dumont \cite{Dumont_79,Dumont_81a} introduced explicitly,
the sequence of polynomials
\be
   X_n(x,y,z)
   \;=\;
   \biggl( yz {\partial \over \partial x}
             \:+\:  xz {\partial \over \partial y}
             \:+\:  xy {\partial \over \partial z} \biggr)^{\! n} \: x
   \;,
 \label{def.schett}
\ee
which we shall call the \textbfit{Schett polynomials} initialized in $x$;
clearly they are homogeneous polynomials of degree $n+1$
with nonnegative integer coefficients,
and are symmetric in $y \leftrightarrow z$.
% The first few Schett polynomials are
% \begin{subeqnarray}
%    X_0  & = &  x   \\
%    X_1  & = &  yz  \\
%    X_2  & = &  x \, (y^2 + z^2)   \\
%    X_3  & = &  yz \, (4 x^2 + y^2 + z^2)   \\
%    X_4  & = &  x \, (4 x^2 y^2 + 4 x^2 z^2 + y^4 + 14 y^2 z^2 + z^4) \\
%    X_5  & = &  yz \, (16 x^4 + 44 x^2 y^2 + 44 x^2 z^2 + y^4 + 14 y^2 z^2 + z^4)
% \end{subeqnarray}
Schett \cite{Schett_76,Schett_77} showed that the polynomials $X_n(x,y,z)$
unify and generalize
the Taylor coefficients of the Jacobian elliptic functions $\sn,\cn,\dn$.
%% see eqns.~\reff{eq.sn.taylor.schett}--\reff{eq.dn.taylor.schett} below.
%% (see Section~\ref{sec.elliptic} for details).

It is easy to see that for even (resp.~odd) $n$,
the polynomial $X_n(x,y,z)$ is divisible by $x$ (resp.~by $yz$)
and that the quotient is a polynomial in $x^2, y^2, z^2$.
We therefore define the \textbfit{reduced Schett polynomials}
\be
   \Xhat_n(x^2,y^2,z^2)
   \;=\;
   \begin{cases}
     \displaystyle {1 \over x} \, X_n(x,y,z)   & \textrm{for $n$ even} \\[4mm]
     \displaystyle {1 \over yz} \, X_n(x,y,z)   & \textrm{for $n$ odd}
   \end{cases}
 \label{def.schett_reduced}
\ee
Then $\Xhat_{2k}$ and $\Xhat_{2k+1}$ are homogeneous polynomials
of degree $k$ in $x^2, y^2, z^2$.
It is also not difficult to show that these two polynomials coincide
when $x=0$:  that is, $\Xhat_{2k}(0,y^2,z^2) = \Xhat_{2k+1}(0,y^2,z^2)$.

\subsection{Total positivity and production matrices}

The main purpose of the present paper is to prove the
{\em coefficientwise Hankel-total positivity}\/
of the even and odd subsequences of Schett polynomials.
Recall first that a finite or infinite matrix of real numbers is called
{\em totally positive}\/ (TP) if all its minors are nonnegative.
% and {\em strictly totally positive}\/ (STP)
% if all its minors are strictly positive.\footnote{
%    {\bf Warning:}  Many authors
%    (e.g.\ \cite{Gantmakher_37,Gantmacher_02,Fomin_00,Fallat_11})
%    use the terms ``totally nonnegative'' and ``totally positive''
%    for what we have termed ``totally positive'' and
%    ``strictly totally positive'', respectively.
%    So it is very important, when seeing any claim about
%    ``totally positive'' matrices, to ascertain~which sense of
%    ``totally positive'' is being used!
%    (This is especially important because many theorems in this subject
%     require {\em strict}\/ total positivity for their validity.)
% }
Background information on totally positive matrices can be found
in \cite{Karlin_68,Gantmacher_02,Pinkus_10,Fallat_11};
they have applications to many areas of pure and applied mathematics.\footnote{
   See \cite[footnote~4]{forests_totalpos} for many references.
%    Including combinatorics
%    \cite{Brenti_89,Brenti_95,Brenti_96,Fomin_00,Skandera_03},
%    stochastic processes \cite{Karlin_59,Karlin_68},
%    statistics \cite{Karlin_68},
%    the mechanics of oscillatory systems \cite{Gantmakher_37,Gantmacher_02},
%    the zeros of polynomials and entire functions
%    \cite{Karlin_68,Asner_70,Kemperman_82,Holtz_03,Pinkus_10,Dyachenko_14},
%    spline interpolation \cite{Schoenberg_53,Karlin_68,Gasca_96},
%    Lie theory \cite{Lusztig_94,Lusztig_98,Fomin_99,Lusztig_08}
%    and cluster algebras \cite{Fomin_10,Fomin_forthcoming},
%    the representation theory of the infinite symmetric group
%    \cite{Thoma_64,Borodin_17},
%    the theory of immanants \cite{Stembridge_91},
%    planar discrete potential theory \cite{Curtis_98,Fomin_01}
%    and the planar Ising model \cite{Lis_17},
%    and several other areas \cite{Gasca_96}.
}
In particular, it is known
\cite[Th\'eor\`eme~9]{Gantmakher_37} \cite[section~4.6]{Pinkus_10}
that an infinite Hankel matrix $(a_{i+j})_{i,j \ge 0}$
of real numbers is totally positive if and only if the underlying sequence
$(a_n)_{n \ge 0}$ is a Stieltjes moment sequence,
i.e.\ the moments of a positive measure on $[0,\infty)$.

However, this is only the beginning of the story,
because our main interest \cite{Sokal_totalpos}
is not with sequences and matrices of real numbers,
but rather with sequences and matrices of {\em polynomials}\/
(with integer or real coefficients) in one or more indeterminates $\bfx$:
in~applications they will typically be generating polynomials that enumerate
some combinatorial objects with respect to one or more statistics.
We equip the polynomial ring $\R[\bfx]$ with the coefficientwise
partial order:  that~is, we say that $P$ is nonnegative
(and write $P \myge 0$)
in~case $P$ is a polynomial with nonnegative coefficients.
%% in the indeterminates~$\bfx$.
We then say that a matrix with entries in $\R[\bfx]$ is
\textbfit{coefficientwise totally positive}
if all its minors are polynomials with nonnegative coefficients;
and we say that a sequence $\ba = (a_n)_{n \ge 0}$ with entries in $\R[\bfx]$
is \textbfit{coefficientwise Hankel-totally positive}
if its associated infinite Hankel matrix
$H_\infty(\ba) = (a_{i+j})_{i,j \ge 0}$
is coefficientwise totally positive.
Most generally, we can consider sequences and matrices
with entries in an arbitrary partially ordered commutative ring;
total positivity and Hankel-total positivity
are then defined in the obvious way.
Coefficientwise Hankel-total positivity of a sequence of polynomials
$(P_n(\bfx))_{n \ge 0}$ {\em implies}\/ the pointwise Hankel-total positivity
(i.e.~the Stieltjes moment property) for all $\bfx \ge 0$,
but it is vastly stronger.

Our first main theorem is:

\begin{theorem}[Coefficientwise Hankel-total positivity of the Schett polynomials]
   \label{thm1.1}
\quad\hfill\vspace*{-1mm}
\begin{itemize}
   \item[(a)]  The sequence $\bigl( \Xhat_{2n}(x^2,y^2,z^2) \bigr)_{n \ge 0}$
      of even reduced Schett polynomials
      is coefficientwise Hankel-totally positive in $x,y,z$.
   \item[(b)]  The sequence $\bigl( \Xhat_{2n+1}(x^2,y^2,z^2) \bigr)_{n \ge 0}$
      of odd reduced Schett polynomials
      is coefficientwise Hankel-totally positive in $x,y,z$.
\end{itemize}
\end{theorem}

\noindent
We remark that the $z=0$ specialization of Theorem~\ref{thm1.1}(a)
follows from a continued fraction due to
Stieltjes \cite[p.~H17]{Stieltjes_1889}
and Rogers \cite[p.~77]{Rogers_07}
together with some general theory \cite{Sokal_totalpos}.\footnote{
   See also Flajolet and Fran\c{c}on \cite[Theorem~1]{Flajolet_89}
   and Milne \cite[Theorems~3.2 and 3.11]{Milne_02}
   for modern presentations of Rogers' elegant proof of the continued fraction.
}
%% \cite[Sections~9.3--9.5]{latpath_SRTR};
%% we will explain this in Section~\ref{??????}.
But the general case of Theorem~\ref{thm1.1} is considerably more subtle.
  
Our proof of Theorem~\ref{thm1.1} will be based on the method of
{\em production matrices}\/ \cite{Deutsch_05,Deutsch_09}.
%% We will review this theory
%% in Sections~\ref{subsec.production}--\ref{subsec.totalpos.prodmat},
%% so now we state only the bare-bones definitions.
Let $P = (p_{ij})_{i,j \ge 0}$ be an infinite matrix
with entries in a commutative ring~$R$;
we assume that $P$ is either row-finite
(i.e.\ has only finitely many nonzero entries in each row)
or column-finite.
Now define an infinite matrix $A = (a_{nk})_{n,k \ge 0}$ by
$a_{nk} = (P^n)_{0k}$.
% \be
%    a_{nk}  \;=\;  (P^n)_{0k}
%    \;.
% \ee
We call $P$ the \textbfit{production matrix}
and $A$ the \textbfit{output matrix}, and we write $A = \scro(P)$.
The two key facts here are the following \cite[section~9.2]{latpath_SRTR}:
%% \cite{Sokal_totalpos}:
if $R$ is a partially ordered commutative ring
and $P$ is totally positive,
then $\scro(P)$ is totally positive
\cite[Theorem~9.4]{latpath_SRTR},
and the zeroth column of $\scro(P)$ is Hankel-totally positive
\cite[Theorem~9.7]{latpath_SRTR}.
%% See Section~\ref{subsec.totalpos.prodmat} for precise statements and proofs.

We will prove Theorem~\ref{thm1.1} by exhibiting production matrices
that generate the even or odd reduced Schett polynomials
as the zeroth column of their output matrix,
and then proving the total positivity of those production matrices.
The production matrices will be
{\em quadridiagonal}\/ unit-lower-Hessenberg matrices:
that is, they will have a superdiagonal with all entries equal to~1,
and then a diagonal, a subdiagonal, and a second subdiagonal;
all other entries are equal to~0.
In the language of \cite{latpath_SRTR},
these are the production matrices for a 2-branched J-fraction.
% {\bf When do they come from a 2-branched S-fraction?????
%    When $y=z$, the even Schett polynomials come from a 2-S-fraction
%    with $\alpha_{3n-1} = (2n)(2n-1)y^2$, $\alpha_{3n} = (2n)^2 x^2$,
%    $\alpha_{3n+1} = (2n)(2n+1)y^2$.
%    And maybe the corresponding odd ones are given by a
%    \emph{modified} 2-S-fraction with the same alphas:
%    compare Proposition~\ref{prop.factorization}.}
Similar but simpler quadridiagonal production matrices
arose in \cite{latpath_laguerre}.
In particular, we will prove the following two results:

\begin{proposition}[Production matrices for even and odd Schett polynomials]
   \label{prop1.2}
\quad\hfill\vspace*{-1mm}
\begin{itemize}
   \item[(a)] The quadridiagonal production matrix
$P(x,y,z) = (p_{ij})_{i,j \ge 0}$ defined by
\begin{subeqnarray}
   p_{n,n+1}  & = &   1   \\[1mm]
   p_{n,n}    & = &   (2n)^2 x^2 \,+\, (2n+1)^2 (y^2 + z^2)  \\[1mm]
   p_{n,n-1}  & = &   (2n)^2 (2n-1) \, \bigl[ (2n-1) x^2 (y^2 + z^2) \,+\,
                                     (2n+1) y^2 z^2 \bigr] \qquad  \\[1mm]
   p_{n,n-2}  & = &   (2n)^2 (2n-2)^2 (2n-1)(2n-3) x^2 y^2 z^2  \\[1mm]
   p_{n,k}    & = &   0 \qquad\textrm{if $k < n-2$ or $k > n+1$}
 \label{eq.prop.prodmat.schetteven}
\end{subeqnarray}
generates the even reduced Schett polynomials
$\bigl( \Xhat_{2n}(x^2,y^2,z^2) \bigr)_{n \ge 0}$
as the zeroth column of its output matrix $\scro(P)$.
   \item[(b)] The quadridiagonal production matrix
$Q(x,y,z) = (q_{ij})_{i,j \ge 0}$ defined by
\begin{subeqnarray}
   q_{n,n+1}  & = &   1   \\[1mm]
   q_{n,n}    & = &   (2n+2)^2 x^2 \,+\, (2n+1)^2 (y^2 + z^2)  \\[1mm]
   q_{n,n-1}  & = &   (2n)^2 (2n+1) \, \bigl[ (2n+1) x^2 (y^2 + z^2) \,+\,
                                     (2n-1) y^2 z^2 \bigr] \qquad \\[1mm]
   q_{n,n-2}  & = &   (2n)^2 (2n-2)^2 (2n+1)(2n-1) x^2 y^2 z^2  \\[1mm]
   q_{n,k}    & = &   0 \qquad\textrm{if $k < n-2$ or $k > n+1$}
 \label{eq.prop.prodmat.schettodd}
\end{subeqnarray}
generates the odd reduced Schett polynomials
$\bigl( \Xhat_{2n+1}(x^2,y^2,z^2) \bigr)_{n \ge 0}$
as the zeroth column of its output matrix $\scro(Q)$.
\end{itemize}
\end{proposition}

\begin{samepage}
\begin{proposition}[Total positivity of the production matrices]
   \label{prop1.3}
\quad\hfill\vspace*{-1mm}
\begin{itemize}
   \item[(a)] The quadridiagonal production matrix $P(x,y,z)$
defined in \reff{eq.prop.prodmat.schetteven}
is coefficientwise totally positive in $x,y,z$.
   \item[(b)] The quadridiagonal production matrix $Q(x,y,z)$
defined in \reff{eq.prop.prodmat.schettodd}
is coefficientwise totally positive in $x,y,z$.
\end{itemize}
\end{proposition}
\end{samepage}

We note that the production matrices $P(x,y,z)$ and $Q(x,y,z)$ coincide
when $x=0$;  this is consistent with the fact that
$\Xhat_{2n}(0,y^2,z^2) = \Xhat_{2n+1}(0,y^2,z^2)$.

\bigskip

Combining Propositions~\ref{prop1.2} and \ref{prop1.3}
with the general theory of totally positive production matrices
\cite[Theorem~9.7]{latpath_SRTR}
proves Theorem~\ref{thm1.1}.

\subsection{Exponential generating functions of the output matrices}

We can also give exponential generating functions for the
output matrices $\scro(P)$ and $\scro(Q)$,
%% arising in Proposition~\ref{prop1.2},
using the theory of exponential Riordan arrays.
%% {\bf Is this only for $\lambda=1$??? Or can we also include $\lambda$???}
Recall first that if $F(t)$ and $G(t)$ are formal power series
with coefficients in a commutative ring $R$, with $G(0) = 0$,
then the \textbfit{exponential Riordan array} $\scrr[F,G]$
is the infinite lower-triangular matrix with entries
\be
   \scrr[F,G]_{nk}
   \;=\;
   {n! \over k!} \:
   [t^n] \, F(t) G(t)^k
   \;\in\; R
   \;.
 \label{def.RFG.0}
\ee
If the ring $R$ contains the rationals (as it will here),
this says that the exponential generating function of
the $k$th column of $\scrr[F,G]$ is $F(t) G(t)^k/k!$.
Note also that if $F$ is even and $G$ is odd,
then $\scrr[F,G]_{nk}$ is nonvanishing
only when $n$ and $k$ have the same parity;
that is, only the even-even and odd-odd submatrices of $\scrr[F,G]$
are nonvanishing.
In this situation we call $\scrr[F,G]$
a \textbfit{checkerboard exponential Riordan array}.

We will show that the output matrix $\scro(P)$ is the even-even submatrix
of a particular checkerboard exponential Riordan array.
Clearly $F$ must be the exponential generating function
of the zeroth column of $\scro(P)$,
which by Proposition~\ref{prop1.2}(a) are the
even reduced Schett polynomials $\Xhat_{2n}(x^2,y^2,z^2)$,
alternated with zero entries.
Hence we must have
\be
   F(t)
   \;=\;
   \sum_{n=0}^\infty \Xhat_{2n}(x^2,y^2,z^2) \: {t^{2n} \over (2n)!}
   \;=\;
   \Dn(t;x,y,z) \, \En(t;x,y,z)
%       \label{def.F.expriordan.P.a} %% \\[3mm]
%    & = &
%    {\cn(t;y',z') \, \dn(t;y',z')
%     \over
%     1 \,-\, x^2 \sn^2(t;y',z')
%    }
   \;,
 \label{def.F.expriordan.P}
\ee
where $\Dn$ and $\En$ are Dumont's hyperelliptic functions
\cite[section~3]{Dumont_81a},
and the final equality is \cite[eq.~(3.11)]{Dumont_81a}.
% (see Section~\ref{subsec.hyperelliptic.dumont}).
% $\sn$, $\cn$ and $\dn$ are the Jacobian elliptic functions
% in Dumont's symmetric parametrization
% (see Section~\ref{subsec.elliptic});
% and $y' = \sqrt{y^2 - x^2}$, $z' = \sqrt{z^2 - x^2}$.
% The second and third equalities in this formula will be explained
% in Section~\ref{subsec.hyperelliptic.dumont}.
For $G$, it turns out that we will have
$G(t) = \Sn(t;x,y,z)$,
% \be
%    G(t)
%    \;=\;
%    \Sn(t;x,y,z)
%       \slabel{def.G.expriordan.P.a} %% \\[2mm]
%    & = &
%    {\sn(t;y',z')
%     \over
%     [1 \,-\, x^2 \sn^2(t;y',z')]^{1/2}
%    }
%    \;,
% \slabel{def.G.expriordan.P.b}
%  \label{def.G.expriordan.P}
% \ee
where $\Sn$ is another of Dumont's hyperelliptic functions.
% and $y'$ and $z'$ are as before;
% the second equality will be explained in
% Section~\ref{subsec.hyperelliptic.dumont}.
We will therefore prove:

\begin{proposition}[Exponential generating function for the output matrix $\scro(P)$]
   \label{prop.expriordan.P}
\hfill\break
Let the matrix $P$ be as in Proposition~\ref{prop1.2}(a),
and let $F(t)$ and $G(t)$ be as above.
%% \reff{def.F.expriordan.P}/\reff{def.G.expriordan.P}.
Then the output matrix $\scro(P)$ is the even-even submatrix
of the exponential Riordan array $\scrr[F,G]$:
\be
   \scro(P)_{nk}
   \;=\;
   {(2n)! \over (2k)!}  \: [t^{2n}] \, F(t) G(t)^{2k}
   \;.
\ee
\end{proposition}

We can do something analogous for the output matrix $\scro(Q)$
associated to the odd Schett polynomials,
but this time it will be most convenient to find it as the {\em odd-odd}\/
submatrix of a checkerboard exponential Riordan array.
We will again have $G(t) = \Sn(t;x,y,z)$.
Since the zeroth column of $\scro(Q)$,
which by Proposition~\ref{prop1.2}(b) are the
odd reduced Schett polynomials $\Xhat_{2n+1}(x^2,y^2,z^2)$,
will now appear in the $k=1$ column of $\scrr[F,G]$
in odd rows starting at $n=1$, we must have
\be
   F(t) G(t)
   \;=\;
   \sum_{n=0}^\infty \Xhat_{2n+1}(x^2,y^2,z^2) \: {t^{2n+1} \over (2n+1)!}
   \;=\;
   \Sn(t;x,y,z) \, \Cn(t;x,y,z)
       \slabel{def.FG.expriordan.Q.a} %% \\[3mm]
%    & = &
%    {\sn(t;y',z')
%     \over
%     1 \,-\, x^2 \sn^2(t;y',z')
%    }
   \;,
% \slabel{def.FG.expriordan.Q.b}
 \label{def.FG.expriordan.Q}
\ee
where $\Sn$ and $\Cn$ are Dumont's hyperelliptic functions,
and the final equality is again \cite[eq.~(3.11)]{Dumont_81a}.
% and $y'$ and $z'$ are as before;
% again the second and third equalities in this formula will be explained
% in Section~\ref{subsec.hyperelliptic.dumont}.
Therefore $F(t) = \Cn(t;x,y,z)$.
% \begin{subeqnarray}
%    F(t)
%    & = &
%    \Cn(t;x,y,z)
%        \slabel{def.F.expriordan.Q.a} \\[2mm]
%    & = &
%    {1
%     \over
%     [1 \,-\, x^2 \sn^2(t;y',z')]^{1/2}
%    }
%    \;.
%  \slabel{def.F.expriordan.Q.b}
%  \label{def.F.expriordan.Q}
% \end{subeqnarray}
We will then prove:

\begin{proposition}[Exponential generating function for the output matrix $\scro(Q)$]
   \label{prop.expriordan.Q}
\hfill\break
Let the matrix $Q$ be as in Proposition~\ref{prop1.2}(b),
and let $F(t)$ and $G(t)$ be as above.
%% \reff{def.G.expriordan.P}/\reff{def.F.expriordan.Q}.
Then the output matrix $\scro(Q)$ is the odd-odd submatrix
of the exponential Riordan array $\scrr[F,G]$:
\be
   \scro(Q)_{nk}
   \;=\;
   {(2n+1)! \over (2k+1)!}  \: [t^{2n+1}] \, F(t) G(t)^{2k+1}
   \;.
\ee
\end{proposition}

\section{Proof of Propositions~\ref{prop1.2}, \ref{prop.expriordan.P} and
   \ref{prop.expriordan.Q}}

We shall prove
Proposition~\ref{prop.expriordan.P} (resp.~\ref{prop.expriordan.Q})
by considering the checkerboard exponential Riordan array defined there
and showing that
the production matrix of its even-even (resp.~odd-odd) submatrix
is indeed what is written in
Proposition~\ref{prop1.2}(a) [resp.~\ref{prop1.2}(b)].
Since the zeroth column of this submatrix is, by construction,
the sequence of even (resp.~odd) reduced Schett polynomials,
this will also prove Proposition~\ref{prop1.2}.

The production matrix of an exponential Riordan array $\scrr[F,G]$
is well known \cite[Theorem~2.19]{forests_totalpos} to be
$p_{nk} = (n!/k!) \: (z_{n-k} \,+\, k \, a_{n-k+1})$,
where $A(s) = \sum_{n=0}^\infty a_n s^n$
and $Z(s) = \sum_{n=0}^\infty z_n s^n$ satisfy
\be
   G'(t) \;=\; A(G(t))  \;,\qquad
   {F'(t) \over F(t)} \;=\; Z(G(t))
   \;.
 \label{eq.prop.riordan.exponential.production.1}
\ee
Here we will give the analogous result
for the even-even and odd-odd submatrices of a
checkerboard exponential Riordan array.

First we define some new series:
% We first express $A(s)$ and $Z(s)$
% in terms of two new series $\Phi(s)$ and $\Psi(s)$,
% as follows \cite{Sokal-Chen_trees_totalpos}:
% \begin{subeqnarray}
%    A(s)   & = &  \Phi(s) \: \Psi(s)
%       \\[1mm]
%    Z(s)   & = &  \Phi(s) \: \Psi'(s)
% \end{subeqnarray}
\be
   B(s)  \;\eqdef\;  A(s)^2   \,,\qquad
   C(s)  \;\eqdef\;  2 A(s) Z(s)  \,,\qquad
   D(s)  \;\eqdef\;  Z(s)^2 + A(s) Z'(s)
   \;.
 \label{def.BCD}
\ee
Since $F$ is even and $G$ is odd,
it follows that $A$ is even and $Z$ is odd;
then $B$ and $D$ are even, while $C$ is odd, and we can write
\be
   B(t)  \;=\;  \sum_{m=0}^\infty b_m \, t^{2m}  \,,\qquad
   C(t)  \;=\;  \sum_{m=0}^\infty c_m \, t^{2m+1}  \,,\qquad
   D(t)  \;=\;  \sum_{m=0}^\infty d_m \, t^{2m}
   \;.
 \label{def.bmcmdm}
\ee
We then have:

\begin{theorem}[Production matrices for submatrices of a checkerboard exponential Riordan array]
   \label{thm.riordan.exponential.production.checkerboard}
Let $\scrr[F,G]$ be a checkerboard exponential Riordan array
with invertible diagonal entries and $F(0) = 1$,
and let $\scrr[F,G]^\eee$ and $\scrr[F,G]^\ooo$
be its even-even and odd-odd submatrices.
Then:
\begin{itemize}
   \item[(a)]  The production matrix of $\scrr[F,G]^\eee$ has entries
\be
   p_{nk}  \;=\;  {(2n)! \over (2k)!} \,
     \Bigl[ 2k(n+k) \, b_{n-k+1} \:+\: 2k \, c_{n-k}  \:+\:  d_{n-k}  \Bigr]
     \;.
 \label{eq.thm.riordan.exponential.production.checkerboard.a}
\ee
   \item[(b)]  The production matrix of $\scrr[F,G]^\ooo$ has entries
\be
   p_{nk}  \;=\;  {(2n+1)! \over (2k+1)!} \,
     \Bigl[ (2k+1)(n+k+1) \, b_{n-k+1} \:+\: (2k+1) \, c_{n-k}  \:+\:  d_{n-k}  \Bigr]
     \;.
 \label{eq.thm.riordan.exponential.production.checkerboard.b}
\ee
\end{itemize}
%% Here $A(s)$ and $Z(s)$ are defined
%% by \reff{eq.prop.riordan.exponential.production.1},
%% and then $B(t), C(t), D(t)$ are defined by \reff{def.BCD}.
\end{theorem}

The proof will be based on the following
simple but powerful result, which we discovered very recently
and which to our surprise is apparently new;
it is inspired by the integration-by-parts arguments employed by
Stieltjes \cite{Stieltjes_1889}, Rogers \cite{Rogers_07}
and others \cite{Flajolet_89,Milne_02}
to deduce classical J-fractions.

\begin{proposition}[Exponential-generating-function method for production matrices]
   \label{prop.prodmat.expgen}
Let $A = (a_{nk})_{n,k \ge 0}$ and $P = (p_{nk})_{n,k \ge 0}$
be row-finite matrices with entries in
a commutative ring $R$ containing the rationals,
with $a_{0k} = \delta_{k0}$.
%% For totalpos.tex, we can let $a_{0k} = u_k$
%    and write $A = \scro(P, {\bf u})$.
Fix integers $r \ge 1$ and $0 \le s < r$,
and define the column exponential generating functions of $A$:
\be
   A_k^{[r,s]}(t)
   \;=\;
   \sum_{n=0}^\infty a_{nk} \, {t^{rn+s} \over (rn+s)!}
   \;.
 \label{eq.prop.prodmat.expgen.1}
\ee
Then the following are equivalent:
\begin{itemize}
   \item[(a)] For all $k \ge 0$,
\be
   {d^r \over dt^r} \, A_k^{[r,s]}(t)
   \;=\;
   \sum_{n=0}^\infty  p_{nk} \, A_n^{[r,s]}(t)
   \;.
 \label{eq.prop.prodmat.expgen.2}
\ee
   \item[(b)] $A = \scro(P)$.
\end{itemize}
\end{proposition}

\proof
On both sides of \reff{eq.prop.prodmat.expgen.2},
the only powers of $t$ that can occur are $t^{rm+s}$ with $m \ge 0$;
for the left-hand side we have here used the hypothesis that $0 \le s < r$.
So take the coefficient of $t^{rm+s}/(rm+s)!$
on both sides of \reff{eq.prop.prodmat.expgen.2}:
\be
   \biggl[ {t^{rm+s} \over (rm+s)!} \biggr]
   \:
   {d^r \over dt^r} \, A_k^{[r,s]}(t)
   \;=\;
   \biggl[ {t^{r(m+1)+s} \over (r(m+1)+s)!} \biggr]
   \:
   A_k^{[r,s]}(t)
   \;=\;
   a_{m+1,k}
\ee
while
\be
   \biggl[ {t^{rm+s} \over (rm+s)!} \biggr]
   \:
   \sum_{n=0}^\infty  p_{nk} \, A_n^{[r,s]}(t)
   \;=\;
   \sum_{n=0}^\infty  p_{nk} \, a_{mn}
   \;.
\ee
But $a_{m+1,k} = \sum\limits_{n=0}^\infty  a_{mn} \, p_{nk}$
is precisely the recurrence stating that $A = \scro(P)$.
\qed

Note that, to prove that $A = \scro(P)$,
it suffices to verify the identities \reff{eq.prop.prodmat.expgen.2}
for one pair $(r,s)$.
On the other hand, if $A = \scro(P)$,
then the identities \reff{eq.prop.prodmat.expgen.2}
hold for all pairs $(r,s)$ satisfying $0 \le s < r$.

We remark also that (a)$\implies$(b) holds even when $s \ge r$;
but in this case the hypothesis \reff{eq.prop.prodmat.expgen.2}
is unlikely to hold, because there are terms that will occur
on the left-hand side (unless they happen to vanish)
that cannot occur on the right-hand side,
namely $t^{rm+s}$ with $m < 0$.
We will now see an example of this.

\sketchofproofof{Theorem~\ref{thm.riordan.exponential.production.checkerboard}}
We apply Proposition~\ref{prop.prodmat.expgen} with $r=2$,
taking $s=0$ for part~(a) and $s=1$ for part~(b).
That is, we will compute the production matrix of the matrix
with elements $a_{nk} = \scrr[F,G]_{2n+s,2k+s}$.
Here we can take in~principle any integer $s \ge 0$,
but we will see at the end why only $s = 0,1$ give a valid final result.

For any $k \in \Z$,
let $H_k(t)$ be the exponential generating function of column $2k+s$
of $\scrr[F,G]$,
\be
   H_k(t)  \;\eqdef\;  {F(t) \, G(t)^{2k+s} \over (2k+s)!}
   \;,
\ee
where by convention we take $H_k(t) \eqdef 0$ if $2k+s < 0$.
A straightforward computation (about 1~page long) using
\reff{eq.prop.riordan.exponential.production.1}--\reff{def.bmcmdm}
shows that
\be
   H''_k(t)  \;=\; \sum_{n=-1}^\infty  p_{nk} \, H_n(t)
 \label{eq.secondproof.expriordan.checkerboard.d}
\ee
where
\be
   p_{nk}
%    & = &
%    {(2n+s)! \over (2k+s)!} \,
%    \Bigl[ \big[ (2k+s)(n-k+1) + (2k+s-1)(2k+s) \big] \, b_{n-k+1}
%        \nonumber \\
%    & & \qquad\qquad\qquad
%            \:+\:
%            (2k+1) c_{n-k} 
%            \:+\:
%            d_{n-k}
%    \Bigr]
%         \\[2mm]
%    & = &
   \;=\;
   {(2n+s)! \over (2k+s)!} \,
   \Bigl[ (2k+s)(n+k+s) \, b_{n-k+1} \:+\: (2k+s) \, c_{n-k}  \:+\:  d_{n-k}
   \Bigr]
   \,.
   \qquad
\ee
In general the term $n=-1$ can contribute to the sum
\reff{eq.secondproof.expriordan.checkerboard.d} when $k=0$,
via $b_{n-k+1} = b_0$
(note that $H_{-1}$ is nontrivial if $s \ge 2$);
but its coefficient is $s(s-1)$ and hence vanishes for $s \in \{0,1\}$
(and in any case $H_{-1} = 0$ for $s \in \{0,1\}$).
Therefore, for $s \in \{0,1\}$ we can apply
Proposition~\ref{prop.prodmat.expgen} with $r=2$.
%% Specializing to $s=0$ gives (a), and specializing to $s=1$ gives (b).
\qed

We also have another proof of
Theorem~\ref{thm.riordan.exponential.production.checkerboard};
we will publish it in the fuller version of this work.

We now apply Theorem~\ref{thm.riordan.exponential.production.checkerboard}
to the checkerboard exponential Riordan arrays defined in
Propositions~\ref{prop.expriordan.P} and \ref{prop.expriordan.Q}.
{}From $G(t) = \Sn(t;x,y,z)$ and the differential equation satisfied by $\Sn$
\cite[eq.~(3.7)]{Dumont_81a},
we have immediately
\be
   A(s)  \;=\;  [(1 + x^2 s^2) (1 + y^2 s^2) (1 + z^2 s^2)]^{1/2}
   \;.
\ee
For Proposition~\ref{prop.expriordan.P} we have
$F(t) = \Dn(t;x,y,z) \, \En(t;x,y,z)$;
using \cite[eq.~(3.3)/(3.4)]{Dumont_81a} we get
\be
   Z(s)  \;=\;  s \, A(s)
     \, \biggl[ {y^2 \over 1 + y^2 s^2} \:+\: {z^2 \over 1 + z^2 s^2} \biggr]
   \;.
\ee
For Proposition~\ref{prop.expriordan.Q} we have
$F(t) = \Cn(t;x,y,z)$;
using \cite[eq.~(3.3)/(3.4)]{Dumont_81a} we get
\be
   Z(s)  \;=\;  s \, A(s) \: {x^2 \over 1 + x^2 s^2}
   \;.
\ee
{}From this we compute, in each case,
the series $B,C,D$ and the coefficients $b_m,c_m,d_m$;
this computation is straightforward and we leave it to the reader.
Substituting these coefficients into
Theorem~\ref{thm.riordan.exponential.production.checkerboard}
yields the production matrices of Proposition~\ref{prop1.2}.

We remark that when $x=y$ and $z=0$,
these two checkerboard exponential Riordan arrays coincide;
this array was found earlier in \cite[eq.~(6.24)]{Deb-Sokal_cycle-alt}.

\section{Proof of Proposition~\ref{prop1.3}}

We begin by observing that the production matrices $P$ and $Q$
defined in Proposition~\ref{prop1.2} have a nice factorization.
Let $T(y,z) \eqdef P(0,y,z) = Q(0,y,z)$ be the production matrix at $x=0$:
it is tridiagonal with matrix elements
\begin{subeqnarray}
   T_{n,n+1}  & = &   1   \\[1mm]
   T_{n,n}    & = &   (2n+1)^2 (y^2 + z^2)  \\[1mm]
   T_{n,n-1}  & = &   4n^2 (4n^2 - 1) y^2 z^2
 \label{def.Tnk}
\end{subeqnarray}
And let $L(x)$ be the lower-bidiagonal matrix with entries
\begin{subeqnarray}
   L_{n,n}  & = &   1   \\[1mm]
   L_{n,n-1}    & = &   4n^2 x^2
 \label{def.L}
\end{subeqnarray}
A simple computation shows:

\enlargethispage{0.3\baselineskip}
\begin{samepage}
\begin{proposition}[Factorization of the production matrices]
   \label{prop.factorization}
\quad\hfill\vspace*{-1mm}
\begin{itemize}
   \item[(a)]  $P(x,y,z) = L(x) \, T(y,z)$.
   \item[(b)]  $Q(x,y,z) = T(y,z) \, L(x)$.
\end{itemize}
\end{proposition}
\end{samepage}
\clearpage

Since a nonnegative bidiagonal matrix is totally positive,
to prove Proposition~\ref{prop1.3} it suffices to show that
$T(y,z)$ is coefficientwise totally positive.
This turns out to be surprisingly difficult;
but it is contained within a tridiagonal special case
of \cite[Theorem~C.1]{latpath_laguerre},
which asserts the coefficientwise total positivity
of a class of quadridiagonal lower-Hessenberg matrices, defined as follows.
Let
\be
   {\bf P}  \;\eqdef\;  L_1 L_2 U \, +\, L_1 D_1 \,+\, L_2 D_2
%            \slabel{eq.thm.prodmat.TP.bis.gen.new3.a} \\[2mm] 
%       & = &    L_1(L_2 U \, +\, D_1) \,+\, L_2 D_2
%            \slabel{eq.thm.prodmat.TP.bis.gen.new3.b} \\[2mm] 
%       & = &    L_2(L_1 U + D_2) \,+\, L_1 D_1 
%            \slabel{eq.thm.prodmat.TP.bis.gen.new3.c}
   \label{eq.thm.prodmat.TP.bis.gen.new3}
\ee
where
\be
   L_1  \;=\;  \alpha I \,+\, \xi L  \,,\qquad
   L_2  \;=\;  \beta I  \,+\, \eta L
   \label{eq.thm.prodmat.TP.bis.gen.new3.L1L2}
\ee
and
\begin{itemize}
	\item $L$ is the lower-bidiagonal matrix with 
		the sequence $a_0,a_1,\ldots$ on the diagonal, 
		the sequence $b_1, b_2, \ldots$ on the subdiagonal,
		and zeroes elsewhere;
	\item $U$ is the upper-bidiagonal matrix with 
		the sequence $c_1, c_2,\ldots$ on the superdiagonal,
		the sequence $d_0, d_1, \ldots$ on the diagonal,
		and zeroes elsewhere;
	\item $D_1$ is the diagonal matrix with entries $e_0, e_1,\ldots\,$;
	\item $D_2$ is the diagonal matrix with entries $f_0, f_1,\ldots\,$;
\end{itemize}
and $\alpha$, $\beta$, $\xi$, $\eta$,
$\bfa = (a_n)_{n\geq 0}$, $\bfb = (b_n)_{n\geq 1}$, 
$\bfc = (c_n)_{n\geq 1}$, $\bfd = (d_n)_{n\geq 0}$,
$\bfe = (e_n)_{n\geq 0}$, $\bfff = (f_n)_{n\ge 0}$
are all indeterminates.
We then have:

\begin{theorem}  {\rm \cite[Theorem~C.1]{latpath_laguerre}}
\label{thm.quadmat2}
The matrix ${\bf P}$ defined by
\reff{eq.thm.prodmat.TP.bis.gen.new3}/%
\reff{eq.thm.prodmat.TP.bis.gen.new3.L1L2}
is totally positive, coefficientwise in the indeterminates
$\alpha,\beta, \xi,\eta, \bfa, \bfb, \bfc, \bfd, \bfe, \bfff$.
\label{thm.quadridiagonal_TP}
\end{theorem}

The quadridiagonal matrices $P$ and $Q$ of Proposition~\ref{prop1.2}
are {\em not}\/ (as far as we can tell) covered by Theorem~\ref{thm.quadmat2},
but the tridiagonal matrix $T$ is.
Taking $\bfd = \bzero$ makes the matrix ${\bf P}$ tridiagonal;
the further specializations
\begin{eqnarray}
   & &
   \alpha = \beta =  1, \quad
   \xi = y^2, \quad
   \eta = z^2,
       \nonumber \\
   & &
   a_n = 0, \quad
   b_n = (2n)(2n+1), \quad
   c_n = 1, \quad
   e_n = (2n+1) z^2, \quad
   f_n = (2n+1) y^2
       \qquad
\end{eqnarray}
then yield $T$.

\sketchofproofof{Theorem~\ref{thm.quadridiagonal_TP}}
We first let
${\bf Q} = \left.{\bf P}\right|_{\bfff = \bzero} = L_1 (L_2 U \, +\, D_1)$.
Next let $\bp_k$, $\bq_k$ and $\bll_k$ denote the $k$-th columns of the
matrices ${\bf P}, {\bf Q}$ and $L_2$, respectively.
%Thus, $\bp_k =  \bq_k +  f_k \bll_k$.

We use induction to prove that for any fixed pair of integers
$0\leq k \leq m+1$,
the matrix $(\bq_0, \bq_1, \ldots, \bq_{k-1}, \bp_{k},\ldots, \bp_m)$
is totally positive.
We do this in the following steps:
\begin{itemize}
   \item Step 1: The matrix ${\bf Q}$ is totally positive:
      this follows from the fact that bidiagonal matrices
      with non-negative entries are totally positive,
      together with the tridiagonal comparison theorem 
      \cite{Sokal_totalpos} \cite[Proposition~3.1]{Zhu_21b}
      \cite[Proposition~2.6]{latpath_laguerre}.
      This establishes the base case $k=m+1$ of our induction.

   \item Step 2: The matrix $(\bq_0, \ldots, \bq_{k-1}, \bll_{k})$
      is totally positive: 
	When $\bfd = \bzero$, this can easily be shown by 
	direct consideration of the minors. 
	The general case is significantly more difficult.

   \item Step 3: If the matrix $(\bp_{k+1}, \ldots, \bp_{m})$
      	is totally positive, then so is \\
      	$(\bq_0, \ldots, \bq_{k-1}, \bll_{k}, \bp_{k+1}, \ldots, \bp_{m})$:
	We prove this by induction where the base case \hbox{$k=m$} 
	is Step 2.
	When $k<m$, we let 
	$\bt_{k+1} =  \bigl(\bp_{k+1}|_{b_{k+1}=0}\bigr)\bigr|_{f_{k+1}\to \xi b_{k+1}c_{k+1}}$
	and let $ (t_{n,k+1})_{n\geq 0} = \bt_{k+1}$.
	Next let $\widetilde{\bp}_{k+1}$ be the same as $\bt_{k+1}$ 
	except that the entry $t_{k,k+1}$ is made equal to $0$. 
	It is not difficult to show that the matrix 
	$(\widetilde{\bp}_{k+1},\bp_{k+2},\ldots,\bp_{m})$ is totally positive.
	Next, we notice that the matrix 
	$S = (\bq_0, \ldots, \bq_{k-1}, \bll_{k}, \widetilde{\bp}_{k+1},\bp_{k+2}, \ldots, \bp_{m})$ 
	consists of two totally positive blocks which overlap in a single row. 
	Using 
	\cite[p.~398]{Karlin_68}, 
	we get that $S$ is totally positive. 
	Finally, we obtain the desired result by right-multiplying $S$ 
	with the upper-bidiagonal matrix that has $1$ on the diagonal, 
	$(\alpha + \xi a_k) c_{k+1}$ in position $(k,k+1)$
	and zeros elsewhere.

   \item Step 4 (Induction step): If the matrix
      $(\bq_0, \bq_1, \ldots, \bq_{k}, \bp_{k+1},\ldots, \bp_m)$
      is totally positive, then so is
      $(\bq_0, \bq_1, \ldots, \bq_{k-1}, \bp_{k},\ldots, \bp_m)$:
	This follows from noticing $\bp_k =  \bq_k +  f_k \bll_k$
	%% to write this matrix as
        %% $(\bq_0, \bq_1, \ldots, \bq_{k-1}, \bq_k +  f_k \bll_k,\bp_{k+1},\ldots, \bp_m)$,
	and then using the column-linearity of determinants.

\end{itemize}
The details of this proof can be found in \cite[Appendix~C]{latpath_laguerre}.
\qed

\section{Remarks}

In a longer version of this work, we intend to treat the following
generalizations:

1) We can define Schett polynomials in any number of variables, not just three.
Fix an integer $m \ge 1$, and define the Schett operator
in $m+1$ variables $x_0,x_1,\ldots,x_m$ by
\be
   \scrd_{m+1}
   \;\eqdef\;
   \sum_{i=0}^m \: \Biggl( \!\! \prod\limits_{\begin{scarray}
                                 0 \le j \le m \\
                                 j \neq i
                              \end{scarray}}
                        \!\!\! x_j \! \Biggr)
                \: {\partial \over \partial x_i}
   \;.
\ee
Then define the
\textbfit{Schett polynomials in $\bm{m+1}$ variables}
initialized in $x_0$ by
\be
   X^{[m+1]}_n(x_0, x_1,\ldots,x_m)
   \;\eqdef\;
   (\scrd_{m+1})^n \, x_0
   \;.
\ee
These polynomials are clearly symmetric in $x_1,\ldots,x_m$.
Since $x_0$ plays a special role, we write
$x = x_0$ and $\bfy = (y_1,\ldots,y_m) = (x_1,\ldots,x_m)$.
%% Specializations of $X^{[m+1]}_n$ to either $x=0$ or $y_i = 0$
%% give the Taylor coefficients of the Jacobian hyperelliptic functions
%% $S^{[m]}$ and $C_i^{[m]}$, respectively
It turns out that if we define $a_i = 1/y_i$
and make an appropriate rescaling,
then the Schett polynomial in $m$ variables
is obtained from the one in $m+1$ variables
by specializing to $a_m = 0$.
This allows us to take $m \to\infty$
and define a symmetric-function generalization of the Schett polynomials.
For this generalization we have obtained analogues of
Propositions~\ref{prop1.2}, \ref{prop.expriordan.P} and \ref{prop.expriordan.Q}
for the case $x=0$,
but the coefficientwise total positivity of this production matrix
is at present a conjecture.
Furthermore, we have not yet been able to find a production matrix
for the case $x \neq 0$.

2) Schett \cite[Theorem~III/1]{Schett_76} showed that $X_n(1,1,1) = n!$.
It was therefore natural to interpret $X_n(x,y,z)$ as enumerating
permutations of $[n] \eqdef \{1,\ldots,n\}$ with respect to some
bivariate statistic, and Dumont \cite{Dumont_79} identified the statistic:
he showed that
\be
   X_n(x,y,z)
   \;=\;
   \begin{cases}
     \displaystyle x \, D_n(x,y,z)   & \textrm{for $n$ even} \\
     \displaystyle y \, D_n(x,y,z)   & \textrm{for $n$ odd}
   \end{cases}
 \label{eq.X.dumont}
\ee
where
\be
   D_n(x,y,z)
   \;\eqdef\;
   \sum_{\sigma \in \Sym_n}
     (x^2)^{\cpeakodd(\sigma)}
     (y^2)^{\cpeakeven(\sigma)}
     z^{\cdrise(\sigma) + \cdfall(\sigma) + \fix(\sigma)}
   \;;
 \label{def.dumontpoly}
\ee
%% --- we call this the \textbfit{Dumont polynomial} ---
here these statistics count, respectively,
odd cycle peaks, even cycle peaks, cycle double rises,
cycle double falls, and fixed points in the permutation $\sigma$.
   (We remark that the $y \leftrightarrow z$ symmetry
   is quite mysterious in this combinatorial interpretation.)
A different combinatorial interpretation of the Schett polynomials,
in terms of vertices of even and odd degree at even and odd levels
of an increasing tree,
was given by Lin and Ma \cite[Theorem~1.6]{Lin_21}.
It will be interesting to find combinatorial interpretations
for the higher-order Schett polynomials.

3) We can generalize \reff{eq.X.dumont}/\reff{def.dumontpoly}
by introducing a factor $\lambda^{\cyc(\sigma)}$,
where $\cyc(\sigma)$ denotes the number of cycles in the permutation $\sigma$.
For these \textbfit{Schett polynomials with cycle-counting},
we have found generalizations of
\reff{eq.prop.prodmat.schetteven}/\reff{eq.prop.prodmat.schettodd}
that empirically generate the correct zeroth-column sequences,
but we have not yet been able to {\em prove}\/ that they do so.
On the other hand, we {\em have}\/ proven the corresponding generalization
of Proposition~\ref{prop1.3}:
these production matrices are totally positive,
coefficientwise in $x,y,z,\lambda$,
as a consequence of the $\bfd = \bzero$ case of Theorem~\ref{thm.quadmat2}.

%% \acknowledgements{ We would like to thank Rod Halburd for sharing with us
%% his knowledge of elliptic and hyperelliptic functions.  ??????????}

%% if you use biblatex then this generates the bibliography
%% if you use some other method then remove this and do it your own way
\enlargethispage{2\baselineskip}
\printbibliography

\end{document}